# MOMENT INEQUALITIES FOR U-STATISTICS[1]

By Radosław Adamczak

*Polish Academy of Sciences*

We present moment inequalities for completely degenerate Banach space valued (generalized) U-statistics of arbitrary order. The estimates involve suprema of empirical processes which, in the real-valued case, can be replaced by simpler norms of the kernel matrix (i.e., norms of some multilinear operators associated with the kernel matrix). As a corollary, we derive tail inequalities for U-statistics with bounded kernels and for some multiple stochastic integrals.

**1. Introduction.** The extensive body of work concerning U-statistics which emerged during the sixty year period following their introduction by Hoeffding has lead to an abundance of results including limit theorems and tail inequalities as well as statistical and combinatorial applications. Most of the results correspond to the classical theorems for sums of independent random variables, exploring the properties of U-statistics under assumptions which are necessary and sufficient for such sums. Although in some cases, such as CLT, those conditions turn out to be necessary and sufficient also for U-statistics, for other problems (like SLLN or LIL) the case of U-statistics is much more complicated and the classical methods of proofs (in particular the existing tail and moment inequalities) are too weak. The properties of U-statistics depend on the so-called *order of degeneracy* and the most troublesome is usually the completely degenerate or canonical case to which other problems can be reduced by means of Hoeffding decomposition (see [10]). It turns out that, already for canonical U-statistics of order 2, what matters is not only the $L^2$- and $L^\infty$-norms of the kernels, but also some more involved norms such as norms of certain operators corresponding to the kernel matrix, as one can see when examining the inequalities by Giné, Latała and Zinn [5]. These quantities have also been reflected in the necessary and sufficient conditions for the LIL for canonical U-statistics of order 2, obtained

Received June 2005; revised May 2006.
[1]Supported in part by MEiN Grant 1 PO3A 012 29.

in [6], and in precise moment estimates for Gaussian chaoses given recently by Latała [9].

In this paper, we generalize the results of [5] to canonical U-statistics of arbitrary order. The organization of the paper is as follows. First, in Section 2, we start from U-statistics with values in a Banach space, then specialize to type 2 spaces. All estimates presented there are expressed in terms of suprema of empirical processes and may be considered counterparts of similar inequalities for Gaussian chaoses due to Borell [3] and Arcones and Giné [2] (see also [1]). The main results are contained in Section 3, where we obtain sharp estimates for moments and tails of canonical U-statistics in the real-valued case. Those estimates involve "deterministic" quantities only and are optimal up to constants and logarithmic factors. Finally, in Section 4, we give analogous tail inequalities for multiple stochastic integrals of bounded deterministic functions with respect to stochastic processes with independent increments and uniformly bounded jumps, in the spirit of inequalities obtained by Houdré and Reynaud-Bouret in [7].

## 2. Estimates involving suprema of empirical processes.

2.1. *Basic definitions and notation.* Let $I_n = \{1, \ldots, n\}$ and consider a measurable space $(\Sigma, \mathcal{F})$ (throughout the paper, we will assume it is a Polish space with the Borel $\sigma$-field) and $(h_{\mathrm{i}})_{\mathrm{i} \in I_n^d}$, a multi-indexed matrix of measurable functions $h_{\mathrm{i}} : \Sigma^d \to B$, for a separable Banach space $(B, |\cdot|)$. Consider also a matrix $(X_i^{(j)})_{i \in I_n, j \in I_d}$ of independent $\Sigma$-valued random variables. To simplify notation, let $h_{\mathrm{i}}$ also stand for $h_{\mathrm{i}}(X_{i_1}^{(1)}, \ldots, X_{i_d}^{(d)})$, where $\mathrm{i} = (i_1, \ldots, i_d)$. Assume that $h_{\mathrm{i}}$ are canonical (completely degenerate), that is, $\mathbb{E}_j h_{\mathrm{i}} = 0$ for all $j \le d$, where $\mathbb{E}_j$ denotes integration with respect to $X^{(j)} = (X_i^{(j)})_{i \in I_n}$ [for $I \subseteq I_d$, we will similarly denote by $\mathbb{E}_I$ integration with respect to $(X_i^{(j)})_{i \in I_n, j \in I}$]. Let us define a random variable

$$Z := \sum_{\mathrm{i} \in I_n^d} h_{\mathrm{i}}(X_{i_1}^{(1)}, X_{i_2}^{(2)}, \ldots, X_{i_d}^{(d)}) = \sum_{\mathrm{i} \in I_n^d} h_{\mathrm{i}}.$$

Our aim is to find precise estimates for the moments of $Z$. To this end, for $J \subseteq I \subseteq I_d$ (not necessarily nonempty) and a fixed value of $\mathrm{i}_{I^c}$, let us introduce the following definition:

DEFINITION 1.

$$\|(h_{\mathrm{i}})_{\mathrm{i}_I}\|_{I,J}$$
$$= \mathbb{E}_{I \setminus J} \sup\Biggl\{ \mathbb{E}_J \sum_{\mathrm{i}_I} \langle \phi, h_{\mathrm{i}} \rangle \prod_{j \in J} f_{i_j}^{(j)}(X_{i_j}^{(j)}) : \phi \in B^*, |\phi| \le 1,$$

(1)

$$\left.\begin{aligned} f^{(j)} = (f_1^{(j)}, \ldots, f_n^{(j)}), f_i^{(j)} : \Sigma \to \mathbb{R}, j \in J, i \in I_n, \\ \mathbb{E}\sum_i |f_i^{(j)}(X_i^{(j)})|^2 \le 1, j \in J \end{aligned}\right\}$$

for $I \neq \varnothing$. Let us further define $\|(h_{\mathrm{i}})_{\mathrm{i}_\varnothing}\|_{\varnothing,\varnothing} = |h_{\mathrm{i}}|$.

Remark. It is worth noting that $\|(h_{\mathrm{i}})_{\mathrm{i}}\|_{I_d,\varnothing} = \mathbb{E}|\sum_{\mathrm{i}} h_{\mathrm{i}}|$. Moreover, for $J \subseteq I = I_d$, $\|(h_{\mathrm{i}})_{\mathrm{i}_I}\|_{I,J}$ is a deterministic quantity (even a norm), whereas for $I \neq I_d$, it is a random variable depending on $(X_{i_k}^{(k)})_{k \in I^c}$.

Remark. Throughout the paper, we use the letter $K$ to denote universal constants, $K_d$ for constants depending on $d$ only and $K_d(B)$ for constants depending on $d$ and some characteristic of a Banach space $B$. In all these cases, the values of constants may differ between occurrences.

2.2. *Inequalities for Banach space valued U-statistics.*

Theorem 1. *There exist constants $K_d$ such that for $p \ge 2$, we have*

$$\mathbb{E}|Z|^p \le K_d^p \left[\sum_{I \subseteq I_d} \sum_{J \subseteq I} p^{p(\#J/2 + \#I^c)} \sum_{\mathrm{i}_{I^c}} \mathbb{E}_{I^c} \|(h_{\mathrm{i}})_{\mathrm{i}_I}\|_{I,J}^p\right]. \tag{2}$$

The main ingredient of the proof of the above theorem is the following lemma which is a corollary of Talagrand's tail inequality for empirical processes [12].

Lemma 1 ([5], Proposition 3.1, see also [4], Theorem 12). *Let $X_1, \ldots, X_n$ be independent random variables with values in $(\Sigma, \mathcal{F})$ and $\mathcal{T}$ a countable class of measurable real functions on $\Sigma$ such that for all $f \in \mathcal{T}$ and $i \in I_n$, $\mathbb{E}f(X_i) = 0$ and $\mathbb{E}f(X_i)^2 < \infty$. Consider the random variable $S := \sup_{f \in \mathcal{T}} |\sum_i f(X_i)|$. Then for all $p \ge 1$,*

$$\mathbb{E}S^p \le K^p\left[(\mathbb{E}S)^p + p^{p/2}\sigma^p + p^p \mathbb{E}\max_i \sup_{f \in \mathcal{T}} |f(X_i)|^p\right],$$

*where*

$$\sigma^2 = \sup_{f \in \mathcal{T}} \sum_i \mathbb{E}f(X_i)^2.$$

To prove Theorem 1, we will need the following simple corollary of Lemma 1.

LEMMA 2. *Let $B$ be a Banach space for which there exists a countable set $D=\{\psi_j\}$ of functionals such that for all $x\in B$,*

$$\|x\|_B = \sup_j |\psi_j(x)|.$$

*Now, let $X_1,\dots,X_n,Y$ be independent random variables with values in $(\Sigma,\mathcal{F})$. Let $E=L^1_Y(B)$, the space of all $B$-valued, integrable functions of the form $f(Y)$ such that $\psi_j\circ f$ is measurable for all $j$. Consider functions $h_i:\Sigma^2\to B$ $(i=1,\dots,n)$ such that $\psi_j\circ h_i$ is measurable for all $j$, $\mathbb{E}_X h_i(X_i,Y)=0$ $Y$-a.e. and $h_i(X_i,Y)\in E$ $X$-a.e. Let $S=\sum_i h_i(X_i,Y)\in E$. Then for all $p\ge 2$,*

$$\begin{aligned}\mathbb{E}\|S\|_E^p &\le K^p\Big[(\mathbb{E}\|S\|_E)^p + p^{p/2}\sigma^p + p^p\mathbb{E}_X\max_i \|h_i(X_i,Y)\|_E^p\Big]\\ &\le K^p\Big[(\mathbb{E}\|S\|_E)^p + p^{p/2}\sigma^p + p^p\mathbb{E}_X\sum_i \|h_i(X_i,Y)\|_E^p\Big],\end{aligned}$$

*where*

$$\begin{aligned}\sigma &= \sup_{f=(f_i(X_i)):\sum\mathbb{E}f_i^2(X_i)\le 1}\mathbb{E}_Y\sup_j\Big|\sum_i \mathbb{E}_X\psi_j(h_i(X_i,Y)f_i(X_i))\Big|\\ &\le \mathbb{E}_Y\sup_j\Big(\sum_i \mathbb{E}_X\psi_j(h_i(X_i,Y))^2\Big)^{1/2}.\end{aligned}$$

PROOF. First, we will construct a countable set of vectors of the form $\phi(Y)=(\phi_1(Y),\phi_2(Y),\dots)$ such that $\sum_j|\phi_j(Y)|=1$ a.e. and for all $g(Y)\in E$,

$$\|g(Y)\|_E = \sup_\phi\Big|\sum_j \mathbb{E}_Y\phi_j(Y)\psi_j(g(Y))\Big|. \tag{3}$$

Note that for every random variable $k(Y)=(k_1(Y),\dots,k_n(Y))\in L^1_Y(\ell^n_\infty)$, there exists a vector $\phi(Y)=(\phi_1(Y),\dots,\phi_n(Y))$ such that

$$\sum_{j=1}^n|\phi_j(Y)|=1 \qquad \text{a.e.} \tag{4}$$

and

$$\mathbb{E}\max_{j\le n}|k_j(Y)| = \Big|\sum_{j=1}^n \mathbb{E}\phi_j(Y)k_j(Y)\Big|,$$

since for each value of $Y$, we can put $\phi_l(Y)=\operatorname{sgn}k_l(Y)$, if $l=\min\{i\le n: |k_i(Y)|=\max_{j\le n}|k_j(Y)|\}$ and $\phi_l(Y)=0$, otherwise. Since all such sequences

$\phi(Y)$, treated as functionals on $L^1_Y(\ell^n_\infty)$, have norm 1 and $L^1_Y(\ell^n_\infty)$ is separable (here, we use the assumption that $\Sigma$ is a Polish space), there exists a countable set $\mathcal{T}_n$ of vectors $\phi(Y)$ satisfying (4) such that

$$\mathbb{E}\max_{j\le n}|k_j(Y)| = \sup_{\phi\in\mathcal{T}_n}\left|\sum_{j=1}^n \mathbb{E}\phi_j(Y)k_j(Y)\right|$$

for all $k(Y)\in L^1_Y(\ell^n_\infty)$. Now, for $g(Y)\in E$,

$$\|g(Y)\|_E = \sup_n \mathbb{E}\max_{j\le n}|\psi_j(g(Y))| = \sup_n\sup_{\phi\in\mathcal{T}_n}\left|\sum_{j\le n}\phi_j(Y)\psi_j(g(Y))\right|,$$

so, to obtain (3), it is enough to take a set consisting of all vectors $\phi(Y)\in\bigcup\mathcal{T}_n$ completed with zeros to vectors of infinite length.

We thus have $S=\sup_\phi|\sum_i\sum_j\mathbb{E}_Y\psi_j(h_i(X_i,Y))\phi_j(Y)| = \sup_\phi|\sum_i g^i_\phi(X_i)|$ and can estimate $\|S\|_p$ using Lemma 1 (although, formally, it deals with the case when the same function is applied to all $X_i$'s, it is easy to see that it also covers our situation). Indeed, we have

$$\begin{aligned}\sigma &= \sup_\phi\left(\sum_i\mathbb{E}_X g^i_\phi(X_i)^2\right)^{1/2}\\ &= \sup_{f=(f_i)\colon\sum_i\mathbb{E}f_i(X_i)^2\le 1}\ \sup_\phi\left|\sum_j\mathbb{E}_Y\sum_i\mathbb{E}_X\psi_j(h_i(X_i,Y)f_i(X_i))\phi_j(Y)\right|\\ &= \sup_f\mathbb{E}_Y\sup_j\left|\sum_i\mathbb{E}_X\psi_j(h_i(X_i,Y)f_i(X_i))\right|.\end{aligned}$$

Proof of Theorem 1. We will proceed by induction with respect to $d$. For $d=1$ the theorem is an obvious corollary of Lemma 1 since

$$\|(h_i)_i\|_{\{1\},\varnothing}=\mathbb{E}Z,$$

$$\|(h_i)_i\|_{\{1\},\{1\}}=\sup_{|\phi|\le 1}\sqrt{\sum_i\mathbb{E}\langle\phi,h_i\rangle^2},$$

$$\sum_i\mathbb{E}\|(h_i)_i\|^p_{\varnothing,\varnothing}\ge\mathbb{E}\max_i|h_i|^p.$$

Let us therefore assume that the inequality is satisfied for all integers smaller then $d$. Let us denote $\tilde{I}^c=I^c\setminus\{d\}$ for $I\subseteq I_d$. The induction assumption for $d_1=d-1$, applied conditionally with respect to $X^{(d)}$, together with Fubini's theorem, implies that

$$\mathbb{E}|Z|^p\le K^p_{d-1}\sum_{I\subseteq\{1,\dots,d-1\}}\sum_{J\subseteq I}\left[p^{p(\#J/2+\#\tilde{I}^c)}\sum_{\mathbf{i}_{\tilde{I}^c}}\mathbb{E}_{\tilde{I}^c}\mathbb{E}_d\left\|\left(\sum_{i_d}h_{\mathbf{i}}\right)_{\mathbf{i}_I}\right\|^p_{I,J}\right].\tag{5}$$

Note that since, in Definition 1, we can restrict the supremum to a countable set of functions, $\mathbb{E}_d\|(\sum_{i_d} h_{\mathrm{i}})_{\mathrm{i}_I}\|_{I,J}^p$ can be estimated by means of Lemma 2 [applied conditionally on $(X_{i_k}^{(k)})_{k\in\tilde{I}^c}$ if $I \neq I_{d-1}$]. We have

$$\begin{aligned} p^{p(\#J/2+\#\tilde{I}^c)} &\sum_{\mathrm{i}_{\tilde{I}^c}} \mathbb{E}_{\tilde{I}^c}\left(\mathbb{E}_d\left\|\left(\sum_{i_d} h_{\mathrm{i}}\right)_{\mathrm{i}_I}\right\|_{I,J}\right)^p \\ &= p^{p(\#J/2+\#(I\cup\{d\})^c)} \sum_{\mathrm{i}_{(I\cup\{d\})^c}} \mathbb{E}_{(I\cup\{d\})^c}\|(h_{\mathrm{i}})_{\mathrm{i}_{I\cup\{d\}}}\|_{I\cup\{d\},J}^p. \end{aligned}$$

Moreover, $\sigma$ from the lemma is bounded by $\|(h_{\mathrm{i}})_{\mathrm{i}_{I\cup\{d\}}}\|_{I\cup\{d\},J\cup\{d\}}$ and

$$\begin{aligned} p^{p(\#J/2+\#\tilde{I}^c)} &\sum_{\mathrm{i}_{\tilde{I}^c}} \mathbb{E}_{\tilde{I}^c} p^{p/2}\|(h_{\mathrm{i}})_{\mathrm{i}_{I\cup\{d\}}}\|_{I\cup\{d\},J\cup\{d\}}^p \\ &= p^{p(\#(J\cup\{d\})/2+\#(I\cup\{d\})^c)} \\ &\quad\times \sum_{\mathrm{i}_{(I\cup\{d\})^c}} \mathbb{E}_{(I\cup\{d\})^c}\|(h_{\mathrm{i}})_{\mathrm{i}_{I\cup\{d\}}}\|_{I\cup\{d\},J\cup\{d\}}^p. \end{aligned}$$

Finally,

$$\begin{aligned} p^{p(\#J/2+\#\tilde{I}^c)} &\sum_{\mathrm{i}_{\tilde{I}^c}} p^p \sum_{i_d} \mathbb{E}_{\tilde{I}^c}\mathbb{E}_d\|(h_{\mathrm{i}})_{\mathrm{i}_I}\|_{I,J}^p \\ &= p^{p(\#J/2+\#I^c)} \sum_{\mathrm{i}_{I^c}} \mathbb{E}_{I^c}\|(h_{\mathrm{i}})_{\mathrm{i}_I}\|_{I,J}^p. \end{aligned}$$

□

THEOREM 2. *Let $B$ be a Banach space of type 2. Then there exist constants $K_d(B)$ depending only on $d$ and the type 2 constant of $B$ such that for all $p \geq 2$,*

$$\begin{aligned} \mathbb{E}|Z|^p \leq K_d(B)^p &\left[\max_{I\subseteq I_d} \mathbb{E}_{I^c} \max_{\mathrm{i}_{I^c}}\left(\sum_{\mathrm{i}_I}\mathbb{E}_I|h_{\mathrm{i}}|^2\right)^{p/2}\right. \\ &\left.+ \sum_{I\subseteq I_d}\sum_{J\subseteq I} p^{p(\#J/2+\#I^c)}\mathbb{E}_{I^c}\max_{\mathrm{i}_{I^c}}\|(h_{\mathrm{i}})_{\mathrm{i}_I}\|_{I,J}^p\right]. \end{aligned} \tag{6}$$

As we can see, the aim is to replace the external sums on the right-hand side of (1) with maxima. To do so, we will use the following lemmas:

LEMMA 3 ([5], inequality (2.6)). *Let $\xi_1,\dots,\xi_N$ be independent, nonnegative random variables. Then for $p > 1$ and $\alpha > 0$, we have*

$$p^{\alpha p}\sum_i \mathbb{E}\xi_i^p \leq 2(1+p^\alpha)\max\left[p^{\alpha p}\mathbb{E}\max_i \xi_i^p, \left(\sum_i \mathbb{E}\xi_i\right)^p\right].$$

LEMMA 4 ([5], Corollary 2.2). *Consider nonnegative kernels* $g_{\mathrm{i}}:\Sigma^d \to \mathbb{R}_+$. *Then for all* $p \geq 1$,

$$\max_{I\subseteq I_d} \mathbb{E}_I \max_{\mathrm{i}_I} \left(\sum_{\mathrm{i}_{I^c}} \mathbb{E}_{I^c} g_{\mathrm{i}}\right)^p \leq \mathbb{E}\left(\sum_{\mathrm{i}\in I_n^d} g_{\mathrm{i}}\right)^p$$

$$\leq K_d^p \sum_{I\subseteq I_d} p^{\#Ip} \mathbb{E}_I \max_{\mathrm{i}_I} \left(\sum_{\mathrm{i}_{I^c}} \mathbb{E}_{I^c} g_{\mathrm{i}}\right)^p.$$

LEMMA 5. *For* $\alpha > 0$, *arbitrary nonnegative kernels* $g_{\mathrm{i}}:\Sigma^d \to \mathbb{R}_+$ *and* $p > 1$, *we have*

$$p^{\alpha p} \sum_{\mathrm{i}\in I_n^d} \mathbb{E} g_{\mathrm{i}}^p \leq K_d^p p^{\alpha d}\left[p^{\alpha p}\mathbb{E}\max_{\mathrm{i}} g_{\mathrm{i}}^p + \sum_{I\subsetneq\{1,\dots,d\}} p^{\#Ip}\mathbb{E}_I \max_{\mathrm{i}_I}\left(\sum_{\mathrm{i}_{I^c}} \mathbb{E}_{I^c} g_{\mathrm{i}}\right)^p\right].$$

PROOF. We use induction with respect to $d$. For $d = 1$, the inequality is implied by Lemma 3. Assume that the lemma is true for all integers smaller than $d$. Applying the induction assumption to $\mathbb{E}_{\{1,\dots,d-1\}}$ and using the same notation as in the proof of Theorem 1, we get

$$p^{\alpha p}\sum_{\mathrm{i}} \mathbb{E} g_{\mathrm{i}}^p \leq K_{d-1}^p p^{\alpha(d-1)} \mathbb{E}_d \sum_{i_d}\Bigg[p^{\alpha p}\mathbb{E}_{\{d\}^c}\max_{\mathrm{i}_{\{d\}^c}} g_{\mathrm{i}}^p$$

$$+ \sum_{I\subsetneq\{1,\dots,d-1\}} p^{\#Ip}\mathbb{E}_I \max_{\mathrm{i}_I}\left(\sum_{\mathrm{i}_{\tilde{I}^c}} \mathbb{E}_{\tilde{I}^c} g_{\mathrm{i}}\right)^p\Bigg].$$

Lemma 3, together with Lemma 4, gives

$$K_{d-1}^p p^{\alpha(d-1)}\mathbb{E}_d \sum_{\mathrm{i}_d} p^{\alpha p}\mathbb{E}_{\{d\}^c}\max_{\mathrm{i}_{\{d\}^c}} g_{\mathrm{i}}^p$$

$$= K_{d-1}^p p^{\alpha(d-1)}\mathbb{E}_{\{d\}^c}\mathbb{E}_d \sum_{i_d} p^{\alpha p}\max_{\mathrm{i}_{\{d\}^c}} g_{\mathrm{i}}^p$$

$$\leq K_d^p p^{\alpha d} p^{\alpha p}\mathbb{E}\max_{\mathrm{i}} g_{\mathrm{i}}^p + K_d^p p^{\alpha d}\mathbb{E}_{\{d\}^c}\left(\sum_{\mathrm{i}} \mathbb{E}_d g_{\mathrm{i}}\right)^p$$

$$\leq K_d^p p^{\alpha d} p^{\alpha p}\mathbb{E}\max_{\mathrm{i}} g_{\mathrm{i}}^p$$

$$+ K_d^p p^{\alpha d}\sum_{I\subseteq\{1,\dots,d-1\}} p^{\#Ip}\mathbb{E}_I\max_{\mathrm{i}_I}\left(\sum_{\mathrm{i}_{I^c}}\mathbb{E}_{I^c} g_{\mathrm{i}}\right)^p.$$

Moreover, for every $I \subsetneq \{1,\dots,d-1\}$, by Lemma 3 (applied with "$\alpha = \#I$"), Lemma 4 and the fact that $p^{\#I} \le p^d \le K_d^p$, we get

$$
\begin{aligned}
K_{d-1}^p p^{\alpha(d-1)} &\sum_{i_d} \mathbb{E}_d p^{\#Ip} \mathbb{E}_I \max_{\mathrm{i}_I} \Biggl(\sum_{\mathrm{i}_{\tilde{I}^c}} \mathbb{E}_{\tilde{I}^c} g_{\mathrm{i}}\Biggr)^p \\
&= K_{d-1}^p p^{\alpha(d-1)} \mathbb{E}_I p^{\#Ip} \sum_{i_d} \mathbb{E}_d \max_{\mathrm{i}_I} \Biggl(\sum_{\mathrm{i}_{\tilde{I}^c}} \mathbb{E}_{\tilde{I}^c} g_{\mathrm{i}}\Biggr)^p \\
&\le K_d^p p^{\alpha(d-1)} p^{\#I} p^{\#Ip} \mathbb{E}_{I\cup\{d\}} \max_{\mathrm{i}_{I\cup\{d\}}} \Biggl(\sum_{\mathrm{i}_{(I\cup\{d\})^c}} \mathbb{E}_{(I\cup\{d\})^c} g_{\mathrm{i}}\Biggr)^p \\
&\quad + K_d^p p^{\alpha(d-1)} p^{\#I} \mathbb{E}_I \Biggl(\sum_{\mathrm{i}} \mathbb{E}_{I^c} g_{\mathrm{i}}\Biggr)^p \\
&\le K_d^p p^{d\alpha} p^{\#Ip} \mathbb{E}_{I\cup\{d\}} \max_{\mathrm{i}_{I\cup\{d\}}} \Biggl(\sum_{\mathrm{i}_{(I\cup\{d\})^c}} \mathbb{E}_{(I\cup\{d\})^c} g_{\mathrm{i}}\Biggr)^p \\
&\quad + K_d^p p^{\alpha d} \sum_{J\subseteq I} p^{\#Jp} \mathbb{E}_J \max_{\mathrm{i}_J} \Biggl(\sum_{\mathrm{i}_{J^c}} \mathbb{E}_{J^c} g_{\mathrm{i}}\Biggr)^p. \qquad \square
\end{aligned}
$$

LEMMA 6. *Let $B$ be a Banach space of type 2. Then there exist constants $K_d(B)$ depending only on $d$ and the type 2 constant of $B$ such that for all $J \subseteq I \subseteq I_d$ and any fixed value of $\mathrm{i}_{I^c}$, one has*

$$\||(h_{\mathrm{i}})_{\mathrm{i}_I}\||_{I,J} \le K_d(B) \sqrt{\sum_{\mathrm{i}_I} \mathbb{E}_I |h_{\mathrm{i}}|^2}.$$

PROOF. The Cauchy–Schwarz inequality gives

$$
\begin{aligned}
\||(h_{\mathrm{i}})_{\mathrm{i}_I}\||_{I,J}^2 &\le \mathbb{E}_{I\setminus J} \sup_{\phi\in B^*, |\phi|\le 1} \mathbb{E}_J \sum_{\mathrm{i}_J} \Biggl\langle \phi, \sum_{\mathrm{i}_{I\setminus J}} h_{\mathrm{i}} \Biggr\rangle^2 \\
&\le \sum_{\mathrm{i}_J} \mathbb{E}_I \Biggl| \sum_{\mathrm{i}_{I\setminus J}} h_{\mathrm{i}} \Biggr|^2 \le K_d(B) \sum_{\mathrm{i}_I} \mathbb{E}_I |h_{\mathrm{i}}|^2. \qquad \square
\end{aligned}
$$

PROOF OF THEOREM 2. One starts from Theorem 1, then applies Lemma 5 for $I \ne I_d$ to $\sum_{\mathrm{i}_{I^c}} \mathbb{E}_{I^c} \||(h_{\mathrm{i}})_{\mathrm{i}_I}\||_{I,J}^p$ with "$p = p/2$" and $\alpha = 2(\#I^c + \#J/2) + \#I^c$ [taking advantage of the fact that $(p/2)^{\alpha d} \le K_d^p$] and finally applies Lemma 6. $\square$

REMARK. From Lemma 4, one can conclude that in cotype 2 spaces (so, in particular, in Hilbert spaces), the quantity $\max_{I\subseteq I_d}\mathbb{E}_{I^c}\max_{\mathrm{i}_{I^c}}(\sum_{\mathrm{i}_I}\mathbb{E}_I|h_{\mathrm{i}}|^2)^{p/2}$ is indispensable (at least up to constants) since one has

$$\mathbb{E}|Z|^p \geq \frac{1}{K_d(B)^p}\mathbb{E}\left(\sum_{\mathrm{i}}|h_{\mathrm{i}}|^2\right)^{p/2} \geq \frac{1}{K_d(B)^p}\max_{I\subseteq I_d}\mathbb{E}_{I^c}\max_{\mathrm{i}_{I^c}}\left(\sum_{\mathrm{i}_I}\mathbb{E}_I|h_{\mathrm{i}}|^2\right)^{p/2}$$

with $K_d(B)$ depending only on $d$ and the cotype 2 constant of $B$.

Let us also consider a corollary of Theorem 2 which is perhaps more "user-friendly." It can easily be obtained by replacing $\|\!|(h_i)_{\mathrm{i}_I}|\!\|_{I,J}$, for $I\neq I_d$ and $\mathbb{E}|Z| = \|\!|(h_{\mathrm{i}})|\!\|_{I_d,\varnothing}$, with the estimates given in Lemma 6.

COROLLARY 1. *If $B$ is of type 2, then there exist constants $K_d(B)$ depending only on $d$ and the type 2 constant of $B$ such that for $p\geq 2$,*

$$\mathbb{E}|Z|^p \leq K_d(B)^p\Bigg[\left(\sum_{\mathrm{i}}\mathbb{E}|h_{\mathrm{i}}|^2\right)^{p/2} + \sum_{J\subseteq I_d, J\neq\varnothing} p^{p\#J/2}\|\!|(h_{\mathrm{i}})|\!\|_{I_d,J}^p$$
$$+ \sum_{I\subsetneq I_d} p^{p(d+\#I^c)/2}\mathbb{E}_{I^c}\max_{\mathrm{i}_{I^c}}\left(\sum_{\mathrm{i}_I}\mathbb{E}_I|h_{\mathrm{i}}|^2\right)^{p/2}\Bigg].$$

**3. The real-valued case.** The purpose of this section is to simplify the estimates of Theorem 2 in the case of real-valued U-statistics. To be more precise, we would like to replace the troublesome suprema of empirical processes $\|(h_{\mathrm{i}})_{\mathrm{i}_I}\|_{I,J}$ by expressions in which the supremum over a class of functions appears outside the expectation. To do so, let us introduce the following definitions:

DEFINITION 2. For a nonempty, finite set $I$, let $\mathcal{P}_I$ be the family consisting of all partitions $\mathcal{J}=\{J_1,\dots,J_k\}$ of $I$ into nonempty, pairwise disjoint subsets. Let us also define for $\mathcal{J}$ (as above), $\deg(\mathcal{J})=k$. Additionally, let $\mathcal{P}_\varnothing=\{\varnothing\}$ with $\deg(\varnothing)=0$.

DEFINITION 3. For a nonempty set $I\subseteq I_d$, consider $\mathcal{J}=\{J_1,\dots,J_k\}\in\mathcal{P}_I$. For an array $(h_{\mathrm{i}})_{\mathrm{i}\in I_n^d}$ of real-valued kernels and any fixed value of $\mathrm{i}_{I^c}$, define

$$\|(h_{\mathrm{i}})_{\mathrm{i}_I}\|_{\mathcal{J}} = \sup\Bigg\{\Bigg|\sum_{\mathrm{i}_I}\mathbb{E}_I h_{\mathrm{i}}(X_{i_1}^{(1)},\dots,X_{i_d}^{(d)})\prod_{j=1}^{\deg(\mathcal{J})} f_{\mathrm{i}_{J_j}}^{(j)}((X_{i_l}^{(l)})_{l\in J_j})\Bigg| :$$
$$\mathbb{E}\sum_{\mathrm{i}_{J_j}}|f_{\mathrm{i}_{J_j}}^{(j)}((X_{i_l}^{(l)})_{l\in J_j})|^2\leq 1 \text{ for } j=1,\dots,\deg(\mathcal{J})\Bigg\}.$$

Moreover, let $\|(h_{\mathrm{i}})_{\mathrm{i}_\varnothing}\|_\varnothing = |h_{\mathrm{i}}|$.

Remark. If $I = I_d$, then the quantity $\|(h_\mathrm{i})_{\mathrm{i}_I}\|_{\mathcal{J}}$ is a deterministic norm, whereas for $I \neq I_d$ it is a random variable depending on $(X_{i_j}^{(j)})_{j\in I^c}$ (one can see that it is just an analogous norm, computed conditionally for a sub-array of smaller dimension).

3.1. *Real U-statistics of order $d = 3$.* First, we will consider the case $d = 3$. Let us adapt the notation to the simplified situation and write

$$Z := \sum_{ijk=1}^{n} h_{ijk}(X_i, Y_j, Z_k),$$

where in all previous definitions, $Y_j, Z_k$ correspond to $X_{i_2}^{(2)}, X_{i_3}^{(3)}$, respectively.

Remark. On closer inspection of Definitions 1 and 3, one can see that

$$|||(h_\mathrm{i})_{\mathrm{i}_I}|||_{I,\varnothing} = \mathbb{E}_I\left|\sum_{\mathrm{i}_I} h_\mathrm{i}\right| \le \sqrt{\sum_{\mathrm{i}_I}\mathbb{E}_I|h_\mathrm{i}|^2} = \|(h_\mathrm{i})_{\mathrm{i}_I}\|_{\{I\}}$$

and for $s \in I$, $|||(h_\mathrm{i})_{\mathrm{i}_I}|||_{I,\{s\}} \le \|(h_\mathrm{i})_{\mathrm{i}_I}\|_{\{I\}}$. Moreover, $|||(h_\mathrm{i})_{\mathrm{i}_I}|||_{I,I} = \|(h_\mathrm{i})_{\mathrm{i}_I}\|_{\mathcal{J}}$, where $\mathcal{J}$ is the partition of $I$ into singletons.

Thus, it follows that to replace all the quantities on the right-hand side of (6) by quantities introduced in Definition 3, one must estimate expressions of the form

$$\begin{aligned}
&|||(h_{ijk})_{ijk}|||_{\{1,2,3\},\{1,2\}}\\
&\quad = \mathbb{E}_Z \sup\left\{\sum_k \sum_{i,j} \mathbb{E}_{X,Y} h_{ijk}(X_i, Y_j, Z_k) f_i(X_i) g_j(Y_j) : \right.\\
&\qquad\qquad \left.\sum_i \mathbb{E} f_i(X_i)^2, \sum_j \mathbb{E} g_j(Y_j)^2 \le 1\right\}.
\end{aligned}$$

Note that one can choose $\mathbf{e}_m = (\mathbf{e}_m^1(X_1), \ldots, \mathbf{e}_m^n(X_n))$, $\mathbf{f}_m = (\mathbf{f}_m^i(Y_i))_{i\le n}$—orthonormal bases in $L^2(X_1) \times \cdots \times L^2(X_n)$ and $L^2(Y_1) \times \cdots \times L^2(Y_n)$, respectively—and denoting $a_{ijk}(Z_k) := \sum_{lm} \mathbb{E}_{XY} h_{lmk}(X_l, Y_m, Z_k)\mathbf{e}_i^l(X_l)\mathbf{f}_j^m(Y_m)$, write

$$\begin{aligned}
|||(h_{ijk})_{ijk}|||_{\{1,2,3\},\{1,2\}} &= \mathbb{E}\sup\left\{\sum_k \sum_{ij} a_{ijk}(Z_k) x_i y_j : \|x\|_2, \|y\|_2 \le 1\right\}\\
&= \mathbb{E}\left\|\sum_k (a_{ijk}(Z_k))_{ij}\right\|_{l_2\to l_2}.
\end{aligned}$$

It turns out that the problem is to estimate the expected operator norm of a sum of independent random matrices (hereafter, we will denote it simply by $\|\cdot\|$, suppressing the index $l_2 \to l_2$).

Before continuing, let us make a few comments concerning the notation. First, to simplify it, we are going to suppress the outer brackets when writing the norms of Definition 3. For example, we will write $\|(h_{ijk})\|_{\{1\}\{2,3\}}$ instead of $\|(h_{ijk})\|_{\{\{1\}\{2,3\}\}}$. Second, note that any array $(a_{ijk}(Z_k))_{ijk}$ corresponds to an array of kernels $(\tilde{h}(g_i^{(1)}, g_j^{(2)}, Z_k))_{ijk} = (a_{ijk}(Z_k) g_i^{(1)} g_j^{(2)})_{ijk}$, where $(g_i^{(1)}, g_j^{(2)})_{ij}$ is an array of independent standard Gaussian random variables. Thus, for any partition $\mathcal{J}$ (as in Definition 3), we can write $\|(a_{ijk}(Z_k))_{\mathrm{i}_I}\|_{\mathcal{J}} = \|(\tilde{h}_{ijk})_{\mathrm{i}_I}\|_{\mathcal{J}}$ (where $g^{(1)}, g^{(2)}$ correspond resp. to $X^{(1)}, X^{(2)}$ of Definition 3). The following proposition explains the connection between these quantities and the corresponding norms of $(h_{ijk})_{ijk}$:

PROPOSITION 1. *For any $\mathcal{J}$, we have $\|(h_{ijk})\|_{\mathcal{J}} = \|(a_{ijk}(Z_k))\|_{\mathcal{J}}$. Moreover,*

$$\|(h_{ijk})\|_{\{1,2,3\}} = \sqrt{\mathbb{E}\sum_{ijk} h_{ijk}^2} = \sqrt{\mathbb{E}\sum_{ijk} a_{ijk}(Z_k)^2},$$

$$\|(h_{ijk})\|_{\{1\}\{2\}\{3\}} = \sqrt{\sup_{\|x\|_2, \|y\|_2 \le 1} \sum_k \mathbb{E}\left(\sum_{ij} a_{ijk}(Z_k) x_i y_j\right)^2},$$

$$\|(h_{ijk})\|_{\{1,3\}\{2\}} = \sqrt{\sup_{\|x\|_2 \le 1} \sum_{ik} \mathbb{E}\left(\sum_{j} a_{ijk}(Z_k) x_j\right)^2},$$

$$\|(h_{ijk})\|_{\{1\}\{2,3\}} = \sqrt{\sup_{\|x\|_2 \le 1} \sum_{jk} \mathbb{E}\left(\sum_{i} a_{ijk}(Z_k) x_i\right)^2}.$$

PROOF. This is a simple fact from theory of $L_2$ spaces, so we will only show the case $\mathcal{J} = \{1,2,3\}$, just to give the flavor of the proof. We have

$$\begin{aligned}\|(h_{\mathrm{i}})\|_{\{1,2,3\}}^2 &= \sum_{ijk} \mathbb{E} h_{ijk}(X_i, Y_j, Z_k)^2 \\ &= \mathbb{E}_Z \sum_k \sup_{\sum_{ij} \mathbb{E} r_{ij}(X_i, Y_j)^2 \le 1} \left| \sum_{ij} \mathbb{E}_{X,Y} h_{ijk}(X_i, Y_j, Z_k) r_{ij}(X_i, Y_j)\right|^2,\end{aligned}$$

but each $(r_{ij})_{i,j=1}^n$ with $\mathbb{E}\sum r_{ij}(X_i, Y_j)^2 \le 1$ can be expressed as

$$r_{ij}(X_i, Y_j) = \sum_{lm} \beta_{lm} e_l^i(X_i) f_m^j(Y_j)$$

with $\sum_{lm}\beta_{lm}^2 \le 1$ [the sum over $l,m$ is, in general, infinite and the equality is satisfied in $\times_{i,j} L_2(X_i,Y_j)$]. Thus,

$$\begin{aligned}
&\sup_{\sum_{ij}\mathbb{E}r_{ij}(X_i,Y_j)^2\le 1}\left|\sum_{ij}\mathbb{E}_{X,Y}h_{ijk}(X_i,Y_j,Z_k)r_{ij}(X_i,Y_j)\right| \\
&\quad = \sup_{\sum_{lm}\beta_{lm}^2\le 1}\left|\sum_{lm}\sum_{ij}\beta_{lm}\mathbb{E}_{X,Y}h_{ijk}(X_i,Y_j,Z_k)e_l^i(X_i)f_m^j(Y_j)\right| \\
&\quad = \sup_{\sum_{lm}\beta_{lm}^2\le 1}\left|\sum_{lm}a_{lmk}(Z_k)\beta_{lm}\right| = \sum_{ij}a_{ijk}(Z_k)^2
\end{aligned}$$

which already implies that

$$\|(h_{ijk})\|_{\{1,2,3\}} = \sqrt{\mathbb{E}\sum_{ijk}a_{ijk}(Z_k)^2}.$$

Thus, it only remains to be shown that $\|(a_{ijk}(Z_k))\|_{\{1,2,3\}}$ also equals the right-hand side of this equality. However, we have

$$\begin{aligned}
&\|(a_{ijk}(Z_k))\|_{\{1,2,3\}}^2 \\
&\quad = \sup_{\sum_{ijk}\mathbb{E}r(g_i^{(1)},g_j^{(2)},Z_k)^2\le 1}\left|\sum_{ijk}\mathbb{E}a_{ijk}(Z_k)g_i^{(1)}g_j^{(2)}r_{ijk}(g_i^{(1)},g_j^{(2)},Z_k)\right| \\
&\quad = \mathbb{E}\sum_{ijk}a_{ijk}(Z_k)^2.
\end{aligned}$$

□

We will also need to introduce an analogue of Definition 3 for deterministic matrices. Note that we can define the norms $\|(a_{ijk})\|_{\mathcal{J}}$ for any deterministic array $(a_{ijk})_{ijk}$ by passing through $(a_{ijk}g_i^1g_j^2g_k^3)_{ijk}$ similarly to the way we did in the case of $(a_{ijk}(Z_k))_{ijk}$. We will, however, follow [9] and give an alternate definition which is equivalent, but more straightforward. Although this section is devoted to U-statistics of order 3, we will consider a more general setting which will also be useful for U-statistics of higher orders.

Definition 4. Let $(a_{\mathbf{i}})_{\mathbf{i}\in I_n^d}$ be a $d$-indexed array of real numbers. For $\mathcal{J}=\{J_1,\ldots,J_k\}\in\mathcal{P}_{I_d}$, we define

$$\|(a_{\mathbf{i}})_{\mathbf{i}}\|_{\mathcal{J}} = \sup\left\{\sum_{\mathbf{i}}a_{\mathbf{i}}x_{\mathbf{i}_{J_1}}^{(1)}\cdots x_{\mathbf{i}_{J_k}}^{(k)} : \sum_{\mathbf{i}_{J_1}}(x_{\mathbf{i}_{J_1}}^{(1)})^2\le 1,\ldots,\sum_{\mathbf{i}_{J_k}}(x_{\mathbf{i}_{J_k}}^{(k)})^2\le 1\right\}.$$

We then have the following:

LEMMA 7 ([9], Theorem 2). *Consider a 3-indexed matrix $A = (a_{ijk})$. Then for any $p \geq 2$,*

$$\mathbb{E}\left\|\left(\sum_k a_{ijk} g_k\right)_{ij}\right\|_{l_2\to l_2} = \mathbb{E}\left\|\left(\sum_k a_{ijk} g_k\right)_{ij}\right\|_{\{1\}\{2\}}$$
$$\leq K\bigg(\|A\|_{\{1\}\{2,3\}} + \|A\|_{\{2\},\{1,3\}}$$
$$+ \frac{1}{\sqrt{p}}\|A\|_{\{1,2,3\}} + \sqrt{p}\|A\|_{\{1\}\{2\}\{3\}}\bigg).$$

REMARK. Although using the same notation for $\|\cdot\|_{\mathcal{J}}$-norms of deterministic arrays and arrays of kernels seems justified by the aforementioned possibility of defining the former via Gaussian chaoses (and also if we interpret them as norms of multilinear operators on proper tensor products of Hilbert spaces), in what follows, we will use Lemma 7 conditionally on the variables $Z_k$. To avoid ambiguity, we will write $\|(a_{ijk}(Z_k))\|_{\mathcal{J},D}$ to stress that we mean a norm of a deterministic array obtained by fixing the random variables $Z_k$.

To proceed, we will need another lemma.

LEMMA 8 ([4], Lemma 7). *Let $X_1,\ldots,X_n$ be independent random variables with values in $(\Sigma,\mathcal{F})$ and let $\mathcal{T}$ be a class of functions $f:\Sigma\to\mathbb{R}$ such that for all $i$, one has $\mathbb{E}f(X_i) = 0$. Then*

$$\mathbb{E}\sup_{f\in\mathcal{T}}\sum_i f^2(X_i) \leq \sup_{f\in\mathcal{T}}\sum_i \mathbb{E}f(X_i)^2 + 32\sqrt{\mathbb{E}M^2}\,\mathbb{E}\sup_{f\in\mathcal{T}}\left|\sum_i f(X_i)\right| + 8\mathbb{E}M^2,$$

*where $M := \max_i \sup_{f\in\mathcal{T}}|f(X_i)|$.*

Now, consider a sequence of independent Rademacher variables $\varepsilon_1,\ldots,\varepsilon_n$, independent of $X, Y, Z$. Using standard symmetrization inequalities, the fact that Rademacher averages are dominated by Gaussian averages and Lemma 7 conditionally on $Z$, we then obtain

$$\mathbb{E}\left\|\left(\sum_k a_{ijk}(Z_k)\right)_{ij}\right\|$$
$$\leq 2\mathbb{E}\left\|\left(\sum_k a_{ijk}(Z_k)\varepsilon_k\right)_{ij}\right\| \leq 2\sqrt{\frac{\pi}{2}}\mathbb{E}\left\|\left(\sum_k a_{ijk}(Z_k) g_k\right)_{ij}\right\|$$
$$\leq K\Big(\mathbb{E}\|(a_{ijk}(Z_k))\|_{\{1\}\{2,3\},D} + \mathbb{E}\|(a_{ijk}(Z_k))\|_{\{2\},\{1,3\},D} \tag{7}$$

$$+\frac{1}{\sqrt{p}}\mathbb{E}\|(a_{ijk}(Z_k))\|_{\{1,2,3\},D}+\sqrt{p}\,\mathbb{E}\|(a_{ijk}(Z_k))\|_{\{1\}\{2\}\{3\},D}\Big).$$

Obviously,

$$\frac{1}{\sqrt{p}}\mathbb{E}\|(a_{ijk}(Z_k))\|_{\{1,2,3\},D}\le\frac{1}{\sqrt{p}}\sqrt{\mathbb{E}\sum_{ijk}a_{ijk}(Z_k)^2},\tag{8}$$

so we are left with the remaining terms.

Let us start with $\mathbb{E}\|(a_{ijk}(Z_k))\|_{\{1\}\{2\}\{3\},D}$. By Lemma 8, we have

$$\begin{aligned}
&\mathbb{E}\|(a_{ijk}(Z_k))\|^2_{\{1\}\{2\}\{3\},D}\\
&\quad=\mathbb{E}\sup_{\|x\|_2,\|y\|_2\le1}\sum_k\Big(\sum_{ij}a_{ijk}(Z_k)x_iy_j\Big)^2\\
&\quad\le\sup_{\|x\|_2,\|y\|_2\le1}\mathbb{E}\sum_k\Big(\sum_{ij}a_{ijk}(Z_k)x_iy_j\Big)^2\\
&\qquad+32\sqrt{\mathbb{E}M^2}\,\mathbb{E}\Big\|\Big(\sum_k a_{ijk}(Z_k)\Big)_{ij}\Big\|+8\mathbb{E}M^2,
\end{aligned}$$

where $M^2=\max_k\sup_{\|x\|_2,\|y\|_2\le1}|\sum_{ij}a_{ijk}(Z_k)x_iy_j|^2=\max_k\|(a_{ijk}(Z_k))_{ij}\|^2$.

We thus obtain

$$\begin{aligned}
&\mathbb{E}\|(a_{ijk}(Z_k))\|_{\{1\}\{2\}\{3\},D}\\
&\quad\le\sqrt{\mathbb{E}\|(a_{ijk}(Z_k))\|^2_{\{1\}\{2\}\{3\},D}}\\
&\quad\le\sqrt{\sup_{\|x\|_2,\|y\|_2\le1}\sum_k\mathbb{E}\Big(\sum_{ij}a_{ijk}(Z_k)x_iy_j\Big)^2}\\
&\qquad+4\sqrt{2}\sqrt{\Big(\mathbb{E}\max_k\|(a_{ijk}(Z_k))_{ij}\|^2\Big)^{1/2}\mathbb{E}\Big\|\Big(\sum_k a_{ijk}(Z_k)\Big)_{ij}\Big\|}\\
&\qquad+2\sqrt{2}\sqrt{\mathbb{E}\max_k\|(a_{ijk}(Z_k))_{ij}\|^2}.
\end{aligned}$$

Now using the inequality $\sqrt{ab}\le\sqrt{p}a/\varepsilon+b\varepsilon/\sqrt{p}$, we finally obtain, for $0<\varepsilon<1$,

$$\begin{aligned}
&\mathbb{E}\|(a_{ijk}(Z_k))\|_{\{1\}\{2\}\{3\},D}\\
&\le K\Bigg(\sqrt{\sup_{\|x\|_2,\|y\|_2\le1}\sum_k\mathbb{E}\Big(\sum_{ij}a_{ijk}(Z_k)x_iy_j\Big)^2}
\end{aligned}\tag{9}$$

$$+ \frac{\varepsilon}{\sqrt{p}}\mathbb{E}\left\|\left(\sum_k a_{ijk}(Z_k)\right)_{ij}\right\| + \frac{\sqrt{p}}{\varepsilon}\sqrt{\mathbb{E}\max_k \|(a_{ijk}(Z_k))_{ij}\|^2}\Bigg).$$

We will now proceed with the term $\mathbb{E}\|(a_{ijk}(Z_k))\|_{\{2\}\{1,3\},D}$.

$$\begin{aligned}
&\mathbb{E}\|(a_{ijk}(Z_k))\|^2_{\{2\}\{1,3\},D} \\
&\quad = \mathbb{E}\sup_{\|y\|_2\le 1}\sum_{i,k}\left(\sum_j a_{ijk}(Z_k)y_j\right)^2 \\
&\quad = \mathbb{E}\sup_{\|y\|_2\le 1}\sum_k\left(\varepsilon_k\sqrt{\sum_i\left(\sum_j a_{ijk}(Z_k)y_j\right)^2}\right)^2 \\
&\quad \le \sup_{\|y\|_2\le 1}\mathbb{E}\sum_{i,k}\left(\sum_j a_{ijk}(Z_k)y_j\right)^2 \\
&\qquad + 32\sqrt{\mathbb{E}\max_k\|(a_{ijk}(Z_k))_{ij}\|^2}\,\mathbb{E}\sup_{\|y\|_2\le 1}\left|\sum_k\varepsilon_k\sqrt{\sum_i\left(\sum_j a_{ijk}(Z_k)y_j\right)^2}\right| \\
&\qquad + 8\mathbb{E}\max_k\|(a_{ijk}(Z_k))_{ij}\|^2,
\end{aligned} \tag{10}$$

where we have again applied Lemma 8, this time to variables $X_k = (Z_k, \varepsilon_k, k)$ and functions $f_y(Z_k, \varepsilon_k, k) = \varepsilon_k\sqrt{\sum_i(\sum_j a_{ijk}(Z_k)y_j)^2}$.

The problem that remains is to estimate the second factor in the product on the right-hand side of the last inequality. Let $g_1, \ldots, g_n$ be independent standard Gaussian random variables, independent of the $Z_k$'s. We have

$$\begin{aligned}
&\mathbb{E}\sup_{\|y\|_2\le 1}\left|\sum_k\varepsilon_k\sqrt{\sum_i\left(\sum_j a_{ijk}(Z_k)y_j\right)^2}\right| \\
&\quad \le \sqrt{\frac{\pi}{2}}\mathbb{E}\sup_{\|y\|_2\le 1}\left|\sum_k g_k\sqrt{\sum_i\left(\sum_j a_{ijk}(Z_k)y_j\right)^2}\right| \\
&\quad = \sqrt{\frac{\pi}{2}}\mathbb{E}\sup_{\|y\|_2\le 1}|X_y| \le 2\sqrt{\frac{\pi}{2}}\mathbb{E}\sup_{\|y\|_2\le 1}X_y,
\end{aligned}$$

where $X_y = \sum_k g_k\sqrt{\sum_i(\sum_j a_{ijk}(Z_k)y_j)^2}$ is a (conditionally) Gaussian process indexed by the $l_2$ unit ball. The covariance structure of $X$ induces a metric on the indexing set, given by

$$d_X(y,\tilde{y})^2 = \mathbb{E}|X_y - X_{\tilde{y}}|^2$$

$$
\begin{aligned}
&= \sum_k \left( \sqrt{\sum_i \Big(\sum_j a_{ijk}(Z_k) y_j\Big)^2} - \sqrt{\sum_i \Big(\sum_j a_{ijk}(Z_k)\tilde{y}_j\Big)^2} \right)^2 \\
&= \sum_k \left( \Big\| \Big(\sum_j a_{ijk}(Z_k) y_j\Big)_i \Big\|_2 - \Big\| \Big(\sum_j a_{ijk}(Z_k)\tilde{y}_j\Big)_i \Big\|_2 \right)^2 \\
&\le \sum_k \Big\| \Big(\sum_j a_{ijk}(Z_k)(y_j - \tilde{y}_j)\Big)_i \Big\|_2^2 \\
&= \sum_{ik} \Big(\sum_j a_{ijk}(Z_k)(y_j - \tilde{y}_j)\Big)^2 = d_{\tilde{X}}(y, \tilde{y})^2,
\end{aligned}
$$

where $\tilde{X}_y = \sum_{ik} g_{ik} \sum_j a_{ijk}(Z_k) y_j$ is another (conditionally) Gaussian process ($g_{ik}$ being i.i.d. standard Gaussian random variables, independent of the $Z_k$'s). Thus, by Slepian's lemma, we get

$$
\begin{aligned}
\mathbb{E} \sup_{\|y\|_2 \le 1} X_y &\le \mathbb{E} \sup_{\|y\|_2 \le 1} \tilde{X}_y \\
&\le \mathbb{E} \sqrt{\sum_j \Big(\sum_{ik} a_{ijk}(Z_k) g_{ik}\Big)^2} \\
&\le \sqrt{\sum_{ijk} \mathbb{E} a_{ijk}(Z_k)^2}.
\end{aligned}
$$

Inserting this inequality into (10) and using the inequality $\sqrt{ab} \le \sqrt{p}a + b/\sqrt{p}$, we eventually obtain

$$
\begin{aligned}
\mathbb{E}\|(a_{ijk}(Z_k))\|_{\{2\}\{1,3\},D} \le{}& \sqrt{\sup_{\|y\|_2 \le 1} \sum_{ik} \mathbb{E}\Big(\sum_j a_{ijk}(Z_k) y_j\Big)^2} \\
&+ \frac{K}{\sqrt{p}} \sqrt{\sum_{ijk} \mathbb{E} a_{ijk}(Z_k)^2} \\
&+ K\sqrt{p} \sqrt{\mathbb{E} \max_k \|(a_{ijk}(Z_k))_{ij}\|^2}.
\end{aligned} \tag{11}
$$

By symmetry,

$$
\mathbb{E}\|(a_{ijk}(Z_k))\|_{\{1\}\{2,3\},D} \le \sqrt{\sup_{\|y\|_2 \le 1} \sum_{jk} \mathbb{E}\Big(\sum_i a_{ijk}(Z_k) y_i\Big)^2}
$$

$$+\frac{K}{\sqrt{p}}\sqrt{\sum_{ijk}\mathbb{E}a_{ijk}(Z_k)^2} \tag{12}$$
$$+K\sqrt{p}\sqrt{\mathbb{E}\max_k\|(a_{ijk}(Z_k))_{ij}\|^2}.$$

Inequalities (7)–(9) (with sufficiently small $\varepsilon$) and (11), (12), together with Proposition 1, yield the following:

THEOREM 3. *For any $p \geq 2$,*

$$\mathbb{E}\left\|\left(\sum_k a_{ijk}(Z_k)\right)_{ij}\right\|_{l_2\to l_2}$$
$$\leq K\Bigg[\frac{1}{\sqrt{p}}\|(a_{ijk}(Z_k))\|_{\{1,2,3\}}+\|(a_{ijk}(Z_k))\|_{\{1\}\{2,3\}}$$
$$+\|(a_{ijk}(Z_k))\|_{\{2\}\{1,3\}}+\sqrt{p}\|(a_{ijk}(Z_k))\|_{\{1\}\{2\}\{3\}}$$
$$+p\sqrt{\mathbb{E}\max_k\|(a_{ijk}(Z_k))_{ij}\|^2}\Bigg].$$

*In particular,*

$$\|(h_{ijk})\|_{\{1,2,3\}\{1,2\}}\leq K\Bigg[\frac{1}{\sqrt{p}}\|(h_{ijk})\|_{\{1,2,3\}}+\|(h_{ijk})\|_{\{1\}\{2,3\}}$$
$$+\|(h_{ijk})\|_{\{2\}\{1,3\}}+\sqrt{p}\|(h_{ijk})\|_{\{1\}\{2\}\{3\}}$$
$$+p\sqrt{\mathbb{E}_Z\max_k\|(h_{ijk})_{ij}\|^2_{\{1\}\{2\}}}\Bigg].$$

Now combining Theorem 2 with Lemma 3 and the remark at the beginning of the present section, we obtain the following theorem:

THEOREM 4. *For any $p \geq 2$, we have*

$$\mathbb{E}\left|\sum_{ijk}h_{ijk}(X_i,Y_j,Z_k)\right|^p$$
$$\leq K^p\left[\sum_{I\subseteq\{1,2,3\}}\sum_{\mathcal{J}\in\mathcal{P}_I}p^{p(\deg(\mathcal{J})/2+\#I^c)}\mathbb{E}_{I^c}\max_{\mathrm{i}_{I^c}}\|(h_{ijk})_{\mathrm{i}_I}\|^p_{\mathcal{J}}\right].$$

3.2. *Real U-statistics of higher order.* To prove a counterpart of Theorem 4, we will need estimates for $|||(h_\mathrm{i})_\mathrm{i}|||_{I_d,I_{d-1}}=\mathbb{E}_d\|(\sum_{i_d}h_\mathrm{i})_{\mathrm{i}_{I_{d-1}}}\|_{\{\{k\}\colon k\in I_{d-1}\}}$. Note that, again, as for $d=3$, by choosing orthonormal bases, we can translate the problem into one of estimating expectation of the norm of a sum

of independent random $(d-1)$-linear operators by the $\|\cdot\|_{\mathcal{J}}$-norms which satisfy a proper version of Proposition 1. The problem thus reduces to estimating $\mathbb{E}\|(\sum_{i_d} a_{\mathrm{i}}(Z_{i_d}))_{\mathrm{i}_{I_{d-1}}}\| = \mathbb{E}\|(\sum_{i_d} a_{\mathrm{i}}(Z_{i_d}))_{\mathrm{i}_{I_{d-1}}}\|_{\{1\},\dots,\{d-1\}}$.

Lemma 9 ([9], Theorem 2). *There exist constants $K_d$ such that for all $p \geq 2$ and any matrix $A = (a_{\mathrm{i}})_{\mathrm{i}\in I_n^d}$,*

$$\mathbb{E}\left\|\left(\sum_{i_d} a_{\mathrm{i}} g_{i_d}\right)_{\mathrm{i}_{I_{d-1}}}\right\|_{\{1\},\dots,\{d-1\}} \leq K_d \sum_{\mathcal{J}\in\mathcal{P}_{I_d}} p^{(1+\deg\mathcal{J}-d)/2}\|(a_{\mathrm{i}})\|_{\mathcal{J}}.$$

Theorem 5. *Let $Z_1,\dots,Z_n$ be independent random variables with values in $(\Sigma,\mathcal{F})$. For $\mathrm{i}\in\mathbb{N}^{d-1}\times I_n$, let $a_{\mathrm{i}}:\Sigma\to\mathbb{R}$ be measurable functions such that $\mathbb{E}_Z a_{\mathrm{i}}(Z_{i_d}) = 0$. There exist constants $K_d$ such that for all $p \geq 2$, we have*

$$\begin{aligned}\mathbb{E}&\left\|\left(\sum_{i_d} a_{\mathrm{i}}(Z_{i_d})\right)_{\mathrm{i}_{I_{d-1}}}\right\| \\ &\leq K_d \sum_{\mathcal{J}\in\mathcal{P}_{I_d}} p^{(1+\deg(\mathcal{J})-d)/2}\|(a_{\mathrm{i}}(Z_{i_d}))_{\mathrm{i}}\|_{\mathcal{J}} \\ &\quad + K_d \sum_{J\in\mathcal{P}_{I_{d-1}}} p^{1+(1+\deg(\mathcal{J})-d)/2}\sqrt{\mathbb{E}\max_{i_d}\|(a_{\mathrm{i}}(Z_{i_d}))_{\mathrm{i}_{I_{d-1}}}\|_{\mathcal{J}}^2},\end{aligned}$$

*where $\|\cdot\|$ denotes the norm of a $(d-1)$-indexed matrix, regarded as a $(d-1)$-linear operator on $(l_2)^{d-1}$ (thus the $\|\cdot\|_{\{1\},\dots,\{d-1\}}$-norm in our notation). In particular,*

$$\begin{aligned}\mathbb{E}&\left\|\left(\sum_{i_d} h_{\mathrm{i}}\right)_{\mathrm{i}_{I_{d-1}}}\right\|_{\{1\},\dots,\{d-1\}} \\ &\leq K_d \sum_{\mathcal{J}\in\mathcal{P}_{I_d}} p^{(1+\deg(\mathcal{J})-d)/2}\|(h_{\mathrm{i}})_{\mathrm{i}}\|_{\mathcal{J}} \\ &\quad + K_d \sum_{J\in\mathcal{P}_{I_{d-1}}} p^{1+(1+\deg(\mathcal{J})-d)/2}\sqrt{\mathbb{E}\max_{i_d}\|(h_{\mathrm{i}})_{\mathrm{i}_{I_{d-1}}}\|_{\mathcal{J}}^2}.\end{aligned}$$

Proof. As in the proof of Theorem 3, we randomize by an independent Rademacher sequence and apply deterministic estimates conditionally on $Z$ (Lemma 9) to obtain

$$\mathbb{E}\left\|\left(\sum_{i_d} a_{\mathrm{i}}(Z_{i_d})\right)_{\mathrm{i}_{I_{d-1}}}\right\| \leq K_d \sum_{\mathcal{J}\in\mathcal{P}_{I_d}} p^{(1+\deg\mathcal{J}-d)/2}\mathbb{E}\|(a_{\mathrm{i}}(Z_{i_d}))\|_{\mathcal{J},D}. \tag{13}$$

Let us consider a general term on the right-hand side of (13), corresponding to $\mathcal{J} = \{J_1, \ldots, J_k\}$ for $\deg(\mathcal{J}) > 1$. Without loss of generality, we can assume that $d \in J_1$. We have (again, by Lemma 8, using arguments similar to those used in the proof of Theorem 3), for $0 < \varepsilon \le 1$,

$$
\begin{aligned}
&p^{(1+k-d)/2}\mathbb{E}\|(a_{\mathrm{i}}(Z_{i_d}))_{\mathrm{i}}\|_{\mathcal{J},D} \\
&\quad \le p^{(1+k-d)/2}\sqrt{\mathbb{E}\sup_{\|(x^{(j)}_{\mathrm{i}_{J_j}})\|_2 \le 1 \,:\, j=2,\ldots,k} \sum_{i_d}\sum_{\mathrm{i}_{J_1\setminus\{d\}}}\left(\sum_{\mathrm{i}_{I_d\setminus J_1}} a_{\mathrm{i}}(Z_{i_d})\prod_{j=2}^{k} x^{(j)}_{\mathrm{i}_{J_j}}\right)^2} \\
&\quad = p^{(1+k-d)/2} \\
&\qquad \times \sqrt{\mathbb{E}\sup_{\|(x^{(j)}_{\mathrm{i}_{J_j}})\|_2 \le 1 \,:\, j=2,\ldots,k} \sum_{i_d}\left\{\varepsilon_{i_d}\left[\sum_{\mathrm{i}_{J_1\setminus\{d\}}}\left(\sum_{\mathrm{i}_{I_d\setminus J_1}} a_{\mathrm{i}}(Z_{i_d})\prod_{j=2}^{k} x^{(j)}_{\mathrm{i}_{J_j}}\right)^2\right]^{1/2}\right\}^2} \\
&\quad \le Kp^{(1+k-d)/2} \\
&\qquad \times \Biggl(\|(a_{\mathrm{i}}(Z_{i_d}))\|_{\mathcal{J}} \\
&\qquad\qquad + \frac{\varepsilon}{\sqrt{p}}\mathbb{E}\sup_{\|(x^{(j)}_{\mathrm{i}_{J_j}})\|_2 \le 1 \,:\, j=2,\ldots,k}\left|\sum_{i_d} g_{i_d}\sqrt{\sum_{\mathrm{i}_{J_1\setminus\{d\}}}\left(\sum_{\mathrm{i}_{I_d\setminus J_1}} a_{\mathrm{i}}(Z_{i_d})\prod_{j=2}^{k} x^{(j)}_{\mathrm{i}_{J_j}}\right)^2}\right| \\
&\qquad\qquad + \frac{\sqrt{p}}{\varepsilon}\sqrt{\mathbb{E}\max_{i_d}\|(a_{\mathrm{i}}(Z_{i_d}))_{\mathrm{i}_{I_{d-1}}}\|^2_{J_1\setminus\{d\},J_2,\ldots,J_k}}\Biggr),
\end{aligned}
$$

where for $J_1 = \{d\}$, we slightly abuse the notation and identify the partition $\{\varnothing, J_2, \ldots, J_k\}$ of $I_{d-1}$ with the partition $\{J_2, \ldots, J_k\}$.

Now, by Slepian's lemma, we obtain (as in the case $d = 3$)

$$
\begin{aligned}
&\mathbb{E}\sup_{\|(x^{(j)}_{\mathrm{i}_{J_j}})\|_2 \le 1 \,:\, j=2,\ldots,k}\left|\sum_{i_d} g_{i_d}\sqrt{\sum_{\mathrm{i}_{J_1\setminus\{d\}}}\left(\sum_{\mathrm{i}_{I_d\setminus J_1}} a_{\mathrm{i}}(Z_{i_d})\prod_{j=2}^{k} x^{(j)}_{\mathrm{i}_{J_j}}\right)^2}\right| \\
&\quad \le 2\mathbb{E}\sup_{\|(x^{(j)}_{\mathrm{i}_{J_j}})\|_2 \le 1 \,:\, j=2,\ldots,k}\left|\sum_{\mathrm{i}_{J_1}} g_{\mathrm{i}_{J_1}}\sum_{\mathrm{i}_{I_d\setminus J_1}} a_{\mathrm{i}}(Z_{i_d})\prod_{j=2}^{k} x^{(j)}_{\mathrm{i}_{J_j}}\right| \\
&\quad \le K_d \sum_{\mathcal{K}\in\mathcal{P}_{I_d},\deg(\mathcal{K})\le k} p^{(1+\deg(\mathcal{K})-k)/2}\mathbb{E}\|(a_{\mathrm{i}}(Z_{i_d}))_{\mathrm{i}}\|_{\mathcal{K},D},
\end{aligned}
$$

where, in the last inequality, we again used Lemma 9.

Note that if $J_1 = \{d\}$, then Slepian's lemma does not change anything (as was the case when $d = 3$), but in order to shorten the (already quite involved) proof, we do not distinguish this case.

Thus, we obtain

$$\begin{aligned} p^{(1+\deg(\mathcal{J})-d)/2}\mathbb{E}\|(a_{\mathrm{i}}(Z_{i_d}))_{\mathrm{i}}\|_{\mathcal{J},D} \\ &\le K_d p^{(1+\deg(\mathcal{J})-d)/2}\|(a_{\mathrm{i}}(Z_{i_d}))_{\mathrm{i}}\|_{\mathcal{J}} \\ &\quad + K_d\varepsilon \sum_{\mathcal{K}\in\mathcal{P}_{I_d},\deg(\mathcal{K})\le k} p^{(1+\deg(\mathcal{K})-d)/2}\mathbb{E}\|(a_{\mathrm{i}}(Z_{i_d}))_{\mathrm{i}}\|_{\mathcal{K},D} \\ &\quad + K_d\varepsilon^{-1}p^{1+(1+\deg(J_1\backslash\{d\},J_2,\ldots,J_k)-d)/2} \\ &\qquad \times \sqrt{\mathbb{E}\max_{i_d}\|(a_{\mathrm{i}}(Z_{i_d}))_{\mathrm{i}_{I_{d-1}}}\|^2_{J_1\backslash\{d\},J_2,\ldots,J_k}}. \end{aligned}$$

The last inequality remains true for $\deg(\mathcal{J}) = 1$ (i.e. for $\mathcal{J} = \{I_d\}$) since $\mathbb{E}\|(a_{\mathrm{i}}(Z_{i_d}))_{\mathrm{i}}\|_{\{I_d\},D} \le \|(a_{\mathrm{i}}(Z_{i_d}))_{\mathrm{i}}\|_{\{I_d\}}$.

Summing over all $\mathcal{J} \in \mathcal{P}_{I_d}$, we get

$$\begin{aligned} \sum_{\mathcal{J}\in\mathcal{P}_{I_d}} p^{(1+\deg(\mathcal{J})-d)/2}\mathbb{E}\|(a_{\mathrm{i}}(Z_{i_d}))_{\mathrm{i}}\|_{\mathcal{J},D} \\ &\le K_d \sum_{\mathcal{J}\in\mathcal{P}_{I_d}} p^{(1+\deg(\mathcal{J})-d)/2}\|(a_{\mathrm{i}}(Z_{i_d}))_{\mathrm{i}}\|_{\mathcal{J}} \\ &\quad + K_d\varepsilon \sum_{\mathcal{J}\in\mathcal{P}_{I_d}} p^{(1+\deg(\mathcal{J})-d)/2}\mathbb{E}\|(a_{\mathrm{i}}(Z_{i_d}))_{\mathrm{i}}\|_{\mathcal{J},D} \\ &\quad + \frac{K_d}{\varepsilon} \sum_{\mathcal{J}\in\mathcal{P}_{I_{d-1}}} p^{1+(1+\deg(\mathcal{J})-d)/2}\sqrt{\mathbb{E}\max_{i_d}\|(a_{\mathrm{i}}(Z_{i_d}))_{\mathrm{i}_{I_{d-1}}}\|^2_{\mathcal{J}}}. \end{aligned}$$

Taking $\varepsilon$ sufficiently small, we obtain a bound for the right-hand side of (13) which allows us to finish the proof. □

Definition 5. We define a partial order $\prec$ on $\mathcal{P}_I$ as

$$\mathcal{I} \prec \mathcal{J}$$

if and only if for all $I \in \mathcal{I}$, there exists $J \in \mathcal{J}$ such that $I \subseteq J$.

Using the basic theory of $L^2$-spaces and Theorem 5, one obtains the following:

Corollary 2. *Let $\mathcal{I} \in \mathcal{P}_{I_{d-1}}$. Then*

$$\mathbb{E}_d\left\|\left(\sum_{i_d} h_{\mathrm{i}}\right)_{\mathrm{i}_{I_{d-1}}}\right\|_{\mathcal{I}}$$

$$\leq \sum_{\mathcal{J}\in\mathcal{P}_{I_d}\,:\,\mathcal{I}\cup\{\{d\}\}\prec\mathcal{J}} p^{(\deg(\mathcal{J})-\deg(\mathcal{I}))/2}\|(h_{\mathrm{i}})_{\mathrm{i}}\|_{\mathcal{J}}$$

$$+ \sum_{\mathcal{J}\in\mathcal{P}_{I_{d-1}}\,:\,\mathcal{I}\prec\mathcal{J}} p^{1+(\deg(\mathcal{J})-\deg(\mathcal{I}))/2}\sqrt{\mathbb{E}_d \max_{i_d}\|(h_{\mathrm{i}})_{\mathrm{i}}\|_{\mathcal{J}}^2}.$$

We would now like to prove Theorem 4 for higher order U-statistics. It turns out that instead of using Theorem 2, it is more convenient to follow its proof and start the induction argument from the very beginning.

LEMMA 10. *There exist constants $K_d$ such that for any $p \geq 2$,*

$$\mathbb{E}\left|\sum_{\mathrm{i}} h_{\mathrm{i}}\right|^p \leq K_d \sum_{I\subseteq I_d}\sum_{\mathcal{J}\in\mathcal{P}_I}\sum_{\mathrm{i}_{I^c}} p^{p(\#I^c+\deg(\mathcal{J})/2)}\mathbb{E}_{I^c}\|(h_{\mathrm{i}})_{\mathrm{i}_I}\|_{\mathcal{J}}^p. \tag{14}$$

PROOF. We employ an easy induction argument in the spirit of the proof of Theorem 1. For $d = 1$, (14) is an immediate consequence of Lemma 1 since $\mathbb{E}|\sum h_i| \leq \sqrt{\mathbb{E}\sum h_i^2} = \|(h_i)_i\|_{\{1\}}$. As for the induction step, one applies the induction assumption (conditionally on $X^{(d)}$) to $\sum_{i_d} h_{\mathrm{i}}$, then uses Lemma 1 and estimates $\mathbb{E}_d\|(\sum_{\mathrm{i}_d} h_{\mathrm{i}})_{\mathrm{i}_I}\|_{\mathcal{J}}$ (for $I \subseteq I_{d-1}, \mathcal{J} \in \mathcal{P}_I$) by means of Corollary 2 [using the fact that $\sqrt{\mathbb{E}_d \max_{i_d}\|(h_{\mathrm{i}})_{\mathrm{i}_{I_{d-1}}}\|_{\mathcal{J}}^2} \leq (\mathbb{E}_d\sum_{i_d}\|(h_{\mathrm{i}})_{\mathrm{i}_{I_{d-1}}}\|_{\mathcal{J}}^p)^{1/p}$]. □

THEOREM 6. *There exist constants $K_d$ such that for $p \geq 2$,*

$$\mathbb{E}\left|\sum_{\mathrm{i}} h_{\mathrm{i}}\right|^p \leq K_d \sum_{I\subseteq I_d}\sum_{\mathcal{J}\in\mathcal{P}_I} p^{p(\#I^c+\deg(\mathcal{J})/2)}\mathbb{E}_{I^c}\max_{\mathrm{i}_{I^c}}\|(h_{\mathrm{i}})_{\mathrm{i}_I}\|_{\mathcal{J}}^p.$$

PROOF. To replace the sums in $\mathrm{i}_{I^c}$ on the right-hand side of (14) with the maximum over $\mathrm{i}_{I^c}$, it is enough to use Lemma 5 for kernels $g_{\mathrm{i}_{I^c}} = \|(h_{\mathrm{i}})_{\mathrm{i}_I}\|_{\mathcal{J}}^2$ with $p/2$ instead of $p$ and $\alpha$ sufficiently large and to notice that for $J \subseteq I^c$ and $\mathcal{J} \in \mathcal{P}_I$, we have $\mathbb{E}_{I^c\setminus J}\sum_{\mathrm{i}_{I^c\setminus J}}\|(h_{\mathrm{i}})_{\mathrm{i}_I}\|_{\mathcal{J}}^2 \leq \|(h_{\mathrm{i}})_{\mathrm{i}_{J^c}}\|_{\{J^c\}}^2$. □

3.3. *Tail estimates for bounded kernels.* Chebyshev's inequality gives the following corollary of Theorem 6:

THEOREM 7. *Assume that all the kernels $h_{\mathrm{i}}$ are bounded. Then there exist constants $K_d$ such that for all $p \geq 2$,*

$$\mathbb{P}\left(\left|\sum_{\mathrm{i}} h_{\mathrm{i}}\right| > K_d\sum_{I\subseteq I_d}\sum_{\mathcal{J}\in\mathcal{P}_I} p^{\#I^c+\deg(\mathcal{J})/2}\max_{\mathrm{i}_{I^c}}\|\|(h_{\mathrm{i}})_{\mathrm{i}_I}\|_{\mathcal{J}}\|_\infty\right) \leq e^{-p}$$

*or, equivalently, for all* $t \geq 0$,

$$\mathbb{P}\left(\left|\sum_{\mathrm{i}} h_{\mathrm{i}}\right| \geq t\right) \leq K_d \exp\left[-\frac{1}{K_d} \min_{I \subseteq I_d, \mathcal{J} \in \mathcal{P}_I} \left(\frac{t}{\| \|(h_{\mathrm{i}})_{\mathrm{i}_I}\|_{\mathcal{J}}\|_\infty}\right)^{2/(\deg(\mathcal{J})+2\#I^c)}\right].$$

Remark. The above theorem is, in a sense, optimal. The recent inequalities for Gaussian chaoses by Latała state that for $h_{\mathrm{i}} = a_{\mathrm{i}} g_{i_1}^{(1)} \cdots g_{i_d}^{(d)}$, we have

$$\mathbb{P}\left(\left|\sum_{\mathrm{i}} h_{\mathrm{i}}\right| \geq k_d \sum_{\mathcal{J} \in \mathcal{P}_{I_d}} p^{\deg(\mathcal{J})/2} \|(h_{\mathrm{i}})_{\mathrm{i}}\|_{\mathcal{J}}\right) \geq k_d \wedge e^{-p},$$

$$\mathbb{P}\left(\left|\sum_{\mathrm{i}} h_{\mathrm{i}}\right| > K_d \sum_{\mathcal{J} \in \mathcal{P}_{I_d}} p^{\deg(\mathcal{J})/2} \|(h_{\mathrm{i}})_{\mathrm{i}}\|_{\mathcal{J}}\right) \leq e^{-p},$$

which shows (together with the CLT for U-statistics) that apart from constants, the components $p^{(\#I^c + \deg(\mathcal{J})/2)} \|(h_{\mathrm{i}})_{\mathrm{i}_I}\|_{\mathcal{J}}$ for $I = I_d$ are correct and cannot be avoided. To discuss the appearance of other components, let us consider a product $V = \mathbb{G} \prod_{i \in I} X_i$, where the $X_i$'s and $\mathbb{G}$ are independent, the $X_i$'s are centered Poisson random variables with parameter 1 and $\mathbb{G} = \sum_{\mathrm{i}_{I^c}} x_{\mathrm{i}_{I^c}} \prod_{j \in I^c} g_{i_j}^{(j)}$ is a Gaussian chaos ($g_i^{(j)}$ are i.i.d. standard Gaussian). Then $V$ is the limit law of U-statistics $V_n$ with kernels $\prod_{i \in I} X_{n,i_j}^{(j)} a_{n,\mathrm{i}_{I^c}} \prod_{j \notin I} g_{i_j}^{(j)}$ ($\mathrm{i} \in I_n^d$), where $X_{n,i_j}^{(j)}$ are centered Bernoulli random variables with parameter $p = 1/n$ and where the coefficients $a_{\mathrm{i}_{I^c}}$ are properly chosen (from the infinite-divisibility of Gaussian variables or by interpreting $\mathbb{G}$ in terms of multiple stochastic integrals). Then $\mathbb{P}(V \geq k_d \alpha_p \sum_{\mathcal{J} \in \mathcal{P}_{I^c}} p^{\deg(\mathcal{J})/2} \|(x_{\mathrm{i}_{I^c}})\|_{\mathcal{J}}) \geq k_d \wedge e^{-p}$, where $\alpha_p^{1/\#I} \log \alpha_p \sim p$, which shows that the other summands are also correct, at least up to a factor of order $(\log p)^{\#I}$.

Further, note that if $X_i^{(j)}$ are i.i.d. random variables and $h_{\mathrm{i}} = h$ for some function $h$, then the quantities appearing in the above theorem simplify, namely $\| \|(h_{\mathrm{i}})_{i_I}\|_{\mathcal{J}}\|_\infty = n^{\#I/2} \| \|h\|_{\mathcal{J}}\|_\infty$. Thus, we obtain the following:

Corollary 3. *If* $h_{\mathrm{i}} = h$ *and* $X_i^{(j)}$ *are i.i.d. random variables, then for any* $t \geq 0$, *we have*

$$\mathbb{P}\left(\left|\sum_{\mathrm{i}} h_{\mathrm{i}}\right| \geq t\right) \leq K_d \exp\left[-\frac{1}{K_d} \min_{I \subseteq I_d, \mathcal{J} \in \mathcal{P}_I} \left(\frac{t}{n^{\#I/2} \| \|h\|_{\mathcal{J}}\|_\infty}\right)^{2/(\deg(\mathcal{J})+2\#I^c)}\right].$$

Remark. In particular, we can see that the tail of the U-statistic generated by a fixed bounded canonical kernel is of order $n^{d/2}$ which agrees with the CLT for such U-statistics. It is also worth pointing out that each of the above theorems has its "undecoupled" version which can be immediately obtained by applying the decoupling results by de la Peña and Montgomery-Smith [11].

**4. Multiple stochastic integrals with respect to stochastic processes with independent increments.** Theorem 6 also yields tail estimates, in the spirit of Theorem 7, for some multiple stochastic integrals (see, e.g., [8] for the necessary definitions). Namely, let $(N_t^{(i)})_{t\in[0,T]}$ $(i \in I_d)$ be independent càdlàg stochastic processes with independent increments, $N_0^{(i)} = 0$. Let $V^i(t) = \operatorname{Var} N_t^{(i)} < \infty$. Moreover, let $\Lambda^i(t) = \mathbb{E}N_t^{(i)}$ be the compensator of $N^{(i)}$ and define $\tilde{N}^{(i)}(t) = N(t) - \Lambda(t)$. Finally, assume that all the jumps of $N^{(i)}$ are uniformly bounded, say by 1, since this is just a matter of normalization and the typical example we have in mind here is the (not necessarily homogeneous) Poisson process.

DEFINITION 6. For a nonempty subset $I \subseteq I_d$ and $\mathcal{J} = \{J_i\}_{i=1}^k \in \mathcal{P}_I$, we define the quantities

$$\|h\|_{\mathcal{J}} = \sup\biggl\{\int_{[0,T]^{\#I}} h(t_1,\dots,t_d) \prod_{j=1}^{\deg(\mathcal{J})} f^{(j)}((t_i)_{i\in J_j}) \prod_{i\in I} dV^i(t_i) : \int_{[0,T]^{\#J_j}} |f^j((t_i)_{i\in J_j})|^2 \prod_{i\in J_j} dV^i(t_i) \le 1\biggr\}.$$

We further define $\|h\|_{\varnothing} = |h|$.

Notice that as in the case of U-statistics, $\|h\|_{\mathcal{J}}$ is a norm when $I = I_d$. Moreover, for $I \neq I_d$, it is a function of $(t_i)_{i\in I^c}$.

We then have the following:

THEOREM 8. *Let $h : [0,T]^d \to \mathbb{R}$ be a bounded Borel measurable function. Consider the stochastic integral*

$$Z = \int_{[0,T]\times\cdots\times[0,T]} h(t_1,\dots,t_d)\, d\tilde{N}_{t_1}^{(1)} \cdots d\tilde{N}_{t_d}^{(d)}.$$

*Then there exist constants $K_d$ such that for all $p \ge 2$,*

$$\mathbb{P}\biggl(|Z| > K_d \sum_{I\subseteq I_d} \sum_{\mathcal{J}\in\mathcal{P}_I} p^{\#I^c + \deg(\mathcal{J})/2} \max_{\mathbf{i}_{I^c}} \|\|h\|_{\mathcal{J}}\|_\infty\biggr) \le e^{-p}.$$

We would like to approximate $h$ by step functions and the stochastic integral by proper U-statistics (or even homogeneous chaoses). However, the best approximation we may hope for is in $L^2$ and almost sure, whereas in Theorems 7 and 8, we have some $L^\infty$-norms. Thus, we must be careful and approximate by step functions $h_n$ for which those norms are bounded by the corresponding norms of $h$. We will use the following:

LEMMA 11. *Consider probability spaces $(\Omega_i,\mu_i)$, $i\le d$, and $\Omega = \times_{i=1,\dots,d}\Omega_i$, $\mu=\bigotimes_{i=1,\dots,d}\mu_i$. Then there exist constants $K_d$ such that for every $\varepsilon>0$ and every measurable subset $A\subseteq\Omega$ with $\mu(A)>1-\varepsilon$, there exists a subset $B\subseteq A$ such that $\mu(B)>1-K_d\varepsilon^{1/2^{d-1}}$ and for all $I\subsetneq I_d$ and $x_I\in\times_{i\in I}\Omega_i$, we have either $B^I_{x_I}=\varnothing$ or $\mu_{I^c}(B^I_{x_I})>1-K_d\varepsilon^{1/2^{d-1}}$, where $\mu_{I^c}=\bigotimes_{i\in I^c}\mu_i$ and $B^I_{x_I}=\{y_{I^c}\in\times_{i\in I^c}\Omega_i : y\in B, y_I=x_I\}$.*

PROOF. Let us first make a comment concerning notation. We will be using induction and, in the process, will be dealing with various subsets $C\subseteq\times_{i\in I}\Omega_i$ for $I\subseteq I_d$. In such a situation, for $J\subsetneq I$ and $x_J\in\times_{i\in J}\Omega_i$, we will denote the set $\{y_{I\setminus J}\in\times_{i\in I\setminus J}\Omega_i : y_I\in C, y_J=x_J\}$ by $C^J_{x_J}$, which may be slightly inconsistent with the notation in the statement of the lemma. Moreover, when writing Cartesian products of several sets, we will pay no attention to the order (regarding the Cartesian product as the set of functions defined on the indexing set and thus making it "commutative").

Let us now proceed with the proof. For $d=1$, the statement is obvious. Let us thus assume that it is true for all numbers smaller than $d>1$. For $\varnothing\ne I\subsetneq I_d$, let $A_I=\{x_I\in\times_{i\in I}\Omega_i : \mu_{I^c}(A^I_{x_I})>1-\sqrt{\varepsilon}\}$. Then by Fubini's theorem, we have $\mu_I(A_I)>1-\sqrt{\varepsilon}$ and by the induction assumption, there exist sets $B_I\subseteq A_I$ with $\mu_I(B_I)>1-K_{d-1}\varepsilon^{1/2^{d-1}}$ and such that all their sections are either empty or of measure greater than $1-K_{d-1}\varepsilon^{1/2^{d-1}}$. Now, let

$$B=\bigcap_{\varnothing\ne I\subsetneq I_d}\bigcup_{z_I\in B_I}\{z_I\}\times A^I_{z_I}=\bigcap_{\varnothing\ne I\subsetneq I_d}\left(B_I\times\left(\times_{i\in I^c}\Omega_i\right)\right)\cap A.$$

We have $\mu(B)>1-K_d\varepsilon^{1/2^{d-1}}$. Let us consider $J\subsetneq I_d$ and arbitrary $x_J$. Then

$$\begin{aligned} B^J_{x_J} &= \{y_{J^c} : y\in B, y_J=x_J\}\\ &= \bigcap_{\varnothing\ne I\subsetneq I_d}\left\{y_{J^c} : y\in\bigcup_{z_I\in B_I}\{z_I\}\times A^I_{z_I}, y_J=x_J\right\}. \end{aligned}\tag{15}$$

We will show that $B^J_{x_J}$ is either empty or of measure greater than $1-K_d\varepsilon^{1/2^{d-1}}$. Assume that there exists $x_{J^c}\in B^J_{x_J}$. Let $x$ be the element of $\times_{i\in I_d}\Omega_i$ given by the "concatenation" of $x_J$ and $x_{J^c}$. Then $x\in B$ and thus $x_I\in B_I$ for all $I\subsetneq I_d$. Thus, for $I\cap J^c\neq\varnothing$, we have $x_{I\cap J^c}\in(B_I)^{I\cap J}_{x_{I\cap J}}$, so this set is nonempty and, as such, by the definition of $B_I$, is of measure greater than $1-K_{d-1}\varepsilon^{1/2^{d-1}}$. Let us now define

$$U=A^J_{x_J}\cap\bigcap_{I\subsetneq I_d, I\cap J^c\neq\varnothing}\left(\bigtimes_{i\in J^c\setminus I}\Omega_i\times(B_I)^{I\cap J}_{x_{I\cap J}}\right).$$

Clearly, $\mu_{J^c}(U)>1-K_d\varepsilon^{1/2^{d-1}}$ since all the intersected sets are of measure greater than $1-K_{d-1}\varepsilon^{1/2^{d-1}}$ (including $A^J_{x_J}$, since $x_J\in B_J\subseteq A_J$). Now, for $x_{J^c}\in U$, we have $x\in A$ (where $x$ is again the "concatenation" of $x_J$ and $x_{J^c}$). Moreover, for $I\cap J^c\neq\varnothing$, $x_{I\cap J^c}\in(B_I)^{I\cap J}_{x_{I\cap J}}$ and thus $x_I\in B_I$. From the discussion following the assumption that $B^J_{x_J}$ is nonempty, this is also the case for $\varnothing\neq I\subseteq J$. Hence, for any $\varnothing\neq I\subsetneq I_d$, we have $x\in\{x_I\}\times A^I_{x_I}$ with $x_I\in B_I$ and so from (15), we have $x_{J^c}\in B^J_{x_J}$. We have thus proved that $U\subseteq B^J_{x_J}$ which implies that $\mu_{J^c}(B^J_{x_J})>1-K_d\varepsilon^{1/2^{d-1}}$. □

LEMMA 12. *There exist step functions, that is, functions of the form*

$$h_n=\sum_{\mathrm{i}\in I^d_{k_n}}a^{(n)}_{\mathrm{i}}\mathbf{1}_{(t^{(n)}_{i_1},t^{(n)}_{i_1+1}]\times\cdots\times(t^{(n)}_{i_d},t^{(n)}_{i_d+1}]},$$

*such that $h_n\to h$ a.e. and in $L^2$ with respect to the product measure on $[0,T]^d$ with marginals determined by $V^i$ and $\|\|h_n\|_{\mathcal J}\|_\infty\le 3\|\|h\|_{\mathcal J}\|_\infty$ for all $I\subsetneq I_d$ and $\mathcal J\in\mathcal P_I$.*

PROOF. First, note that if we replace $N^{(i)}$ with $c_iN^{(i)}$, then $\|h\|_{\mathcal J}$ multiplies by $\prod_{i\in I}c_i$, so without loss of generality, we can assume that $V^{(i)}(T)=1$ which will allow us to use Lemma 11. Consider any sequence $\tilde h_n$ of step functions converging a.e. to $h$ with $\|\tilde h_n\|_\infty\le\|h\|_\infty$. For any $I\subsetneq I_d$ and $\mathcal J\in\mathcal P_I$, we have $\|\tilde h_n\|_{\mathcal J}\to\|h\|_{\mathcal J}$ a.e., thus we can pass to a subsequence and assume that for a large subset $A^{(n)}_{I^c}$ [say $(A^{(n)}_{I^c})^c$ with measure smaller than $\varepsilon/2^{n2^{d-1}}$, $\varepsilon$ to be chosen later], we have $\|\tilde h_n\|_{\mathcal J}\le 2\|\|h\|_{\mathcal J}\|_\infty$. Then let $B^{(n)}$ be the subset of $\bigcap_I(A^{(n)}_{I^c}\times[0,T]^I)$ given by Lemma 11 applied to the $\sigma$-field generated by sets of the form $(t^{(n)}_{i_1},t^{(n)}_{i_1+1}]\times\cdots\times(t^{(n)}_{i_d},t^{(n)}_{i_d+1}]$, where the $t^{(n)}_i$ correspond to the step function $h_n$, as in the formulation of the lemma. Define for $t=(t_1,\dots,t_d)$, $h_n(t)=\tilde h_n(t)\mathbf{1}_{B^{(n)}}(t)$. Then $h_n$ is a step function and by the Borel–Cantelli lemma, we still have $h_n\to h$ a.e. and in $L^2$ (by the Lebesgue dominated convergence theorem). Moreover, for all $I$ and $t_{I^c}$, the function $g_n(t_I)=h_n(t)$ either equals 0 or differs from $\tilde h_n$ on the set of measure not greater then $K_d\varepsilon^{1/2^{d-1}}/2^n$. If $g_n$ does not equal 0, we have $t_{I^c}\in A^{(n)}_{I^c}$ and thus $\|\tilde{h}_n\|_{\mathcal{J}} \le 2\|\|h\|_{\mathcal{J}}\|_\infty$ at $t_{I^c}$. Thus, $\|h_n\|_{\mathcal{J}} = 0$ or $\|h_n\|_{\mathcal{J}} \le \|h_n - \tilde{h}_n\|_{\mathcal{J}} + \|\tilde{h}_n\|_{\mathcal{J}} \le K_d(\|h\|_\infty + \|\tilde{h}_n\|_\infty)\varepsilon^{1/2^d}/2^{n/2} + 2\|\|h\|_{\mathcal{J}}\|_\infty \le 3\|\|h\|_{\mathcal{J}}\|_\infty$ for $\varepsilon$ sufficiently small. □

Proof of Theorem 8. We will prove moment inequalities for stochastic integrals of bounded kernels. These will imply the theorem by the Chebyshev inequality. Consider functions $h_n$, given by Lemma 12. We can assume that $\max_{i\le k_n} |t_{i+1}^{(n)} - t_i^{(n)}| \to_n 0$. Let $Z_n$ be the $d$-fold stochastic integral of $h_n$. Since $h_n \to h$ in $L^2$, we have $Z_n \to Z$ in $L^2$ and we can assume that $Z_n \to Z$ a.e. Let us now (with a slight abuse of notation) denote by $\|Z_n\|_{\mathcal{J}}$ the $\|\cdot\|_{\mathcal{J}}$-norms of the matrix of kernels which define the homogeneous chaos $Z_n$ viewed as a U-statistic (to distinguish them from $\|h_n\|_{\mathcal{J}}$ given in Definition 6). One can see that for $\mathcal{J} \in \mathcal{P}_{I_d}$, we have $\|Z_n\|_{\mathcal{J}} \le \|h_n\|_{\mathcal{J}}$ and for $I \subsetneq I_d$, $\mathcal{J} \in \mathcal{P}_I$, any fixed value of $\mathrm{i}_{I^c}$ and each $t_{I^c} \in \times_{k\in I^c}(t_{i_k}^{(n)}, t_{i_k+1}^{(n)}]$, we have

$$\|Z_n\|_{\mathcal{J}} \le \|h_n\|_{\mathcal{J}} \prod_{k\in I^c} \left|\tilde{N}^{(k)}_{t^{(n)}_{i_k+1}} - \tilde{N}^{(k)}_{t^{(n)}_{i_k}}\right|,$$

where $\|h_n\|_{\mathcal{J}}$ on the right-hand side is taken at the point $t_{I^c}$. Thus, by Fatou's lemma, Theorem 6 and the definition of $h_n$, we get

$$\begin{aligned}
\mathbb{E}|Z|^p &= \mathbb{E}\liminf_n |Z_n|^p \le \liminf_n \mathbb{E}|Z_n|^p \\
&\le \liminf_n K_d^p \Biggl( \sum_{\mathcal{J}\in\mathcal{P}_{I_d}} p^{p\deg(\mathcal{J})/2} \|h_n\|_{\mathcal{J}}^p \\
&\qquad + \sum_{I\subsetneq I_d} \sum_{\mathcal{J}\in\mathcal{P}_I} p^{p(\#I^c + \deg(\mathcal{J})/2)} \\
&\qquad\qquad \times \mathbb{E}_{I^c} \max_{\mathrm{i}_{I^c}} \|\|h_n\|_{\mathcal{J}}\|_\infty^p \prod_{k\in I^c} \left|\tilde{N}^{(k)}_{t^{(n)}_{i_k+1}} - \tilde{N}^{(k)}_{t^{(n)}_{i_k}}\right|^p \Biggr) \\
&\le K_d^p \sum_{I\subseteq I_d} \sum_{\mathcal{J}\in\mathcal{P}_I} p^{p(\#I^c+\deg(\mathcal{J})/2)} \|\|h\|_{\mathcal{J}}\|_\infty^p \\
&\qquad \times \mathbb{E}_{I^c} \limsup_n \max_{\mathrm{i}_{I^c}} \prod_{k\in I^c} \left|\tilde{N}^{(k)}_{t^{(n)}_{i_k+1}} - \tilde{N}^{(k)}_{t^{(n)}_{i_k}}\right|^p \\
&\le K_d^p \sum_{I\subseteq I_d} \sum_{\mathcal{J}\in\mathcal{P}_I} p^{\#I^c+\deg(\mathcal{J})/2} \|\|h\|_{\mathcal{J}}\|_\infty^p,
\end{aligned}$$

where, in the two last inequalities, we have used the assumption that the jumps of $N^{(k)}$ are bounded by 1 (since the lim sup is then also bounded by a constant and, moreover, the processes $N^{(k)}$ have all moments, which together with Doob's inequality allows us to use Fatou's lemma for lim sup). □

**Acknowledgment.** The author would like to thank Prof. Rafał Latała for all the valuable remarks which influenced the final results of this article.

INSTITUTE OF MATHEMATICS
POLISH ACADEMY OF SCIENCES
ŚNIADECKICH 8
P.O. BOX 21
00-956 WARSZAWA
POLAND
E-MAIL: R.Adamczak@impan.gov.pl